\documentclass[11pt]{amsart}

\usepackage{amsmath,amssymb,amsthm}
\usepackage{xcolor}
\usepackage[margin=1.15in]{geometry}
\usepackage[colorlinks=true,citecolor=blue,linkcolor=blue,urlcolor=blue]{hyperref}

\newtheorem{theorem}{Theorem}[section]
\newtheorem{conjecture}[theorem]{Conjecture}
\newtheorem{proposition}[theorem]{Proposition}
\newtheorem{lemma}[theorem]{Lemma}
\newtheorem{corollary}[theorem]{Corollary}
\theoremstyle{definition}
\newtheorem{definition}[theorem]{Definition}

\newcommand{\Kbar}{\overline K}

\newcommand{\Sp}{\operatorname{Sp}}
\newcommand{\Span}{\operatorname{Sp}_{\mathrm{an}}}
\newcommand{\CC}{\mathbb C}

\title[Geometric Bombieri--Lang]{Recent progress on the geometric Bombieri--Lang conjecture}
\author{Junyi Xie}
\address{Beijing International Center for Mathematical Research, Peking University, Beijing 100871, China}
\email{xiejunyi@bicmr.pku.edu.cn}
\date{}

\begin{document}

\begin{abstract}
We survey recent progress on the geometric Bombieri--Lang conjecture over
function fields of characteristic zero.  We discuss recent work of Xie--Yuan
and Guoquan Gao, which together proves the conjecture for varieties admitting
finite morphisms to abelian varieties.  The guiding idea, developed in joint
work with Xinyi Yuan, is that Vojta's dictionary can be made concrete in this
setting: from rational points of large height one constructs entire curves on
complex fibers.
\end{abstract}

\maketitle
\setcounter{tocdepth}{1}
\tableofcontents

\section{Introduction}

\subsection{Mordell's conjecture and Bombieri--Lang}

One of the starting points of modern Diophantine geometry is the Mordell
conjecture.  It asserts that a smooth projective curve of genus at least two
over a number field has only finitely many rational points.  This theorem,
proved by Faltings \cite{FaltingsMordell}, was a milestone in the subject.  It
was later reproved by Vojta via Diophantine approximation
\cite{VojtaMordell}, and by Lawrence--Venkatesh using $p$-adic Hodge theory
\cite{LawrenceVenkatesh}.

The Bombieri--Lang conjecture is a higher-dimensional continuation of this
story.  A smooth projective curve has genus at least two precisely when it is
of general type.  Thus the most direct higher-dimensional analogue predicts
the following statement.

\begin{conjecture}[Bombieri--Lang, non-density form]\label{conj:BL-density-intro}
Let $X$ be a projective variety of general type over a number field $K$.  Then
$X(K)$ is not Zariski dense in $X$.
\end{conjecture}

This is a qualitative statement: it says that rational points are constrained,
but it does not describe exactly where they should lie.

A more structural form of the conjecture uses the special locus.  Very
roughly, the special locus is the part of a variety swept out by images of
abelian varieties.  Such images are expected sources of many rational points
and also of entire curves.  Therefore one should not expect a finiteness
statement along them.  The special-locus expectation is that rational
points should be finite outside the special locus.  This more explicit form
will be recalled in Section~\ref{sec:conjectures}.

\subsection{New phenomena over function fields}

Now let $K/k$ be a finitely generated extension of characteristic zero.  One
may think first of the case where $K=k(B)$, with $B$ a smooth projective curve
over $k$.  A new phenomenon immediately appears.  If
$X=X_0\times_k K$ is a constant variety, then every point of $X_0(k)$ gives a
$K$-rational point of $X$.  Thus one cannot expect a literal finiteness
statement for $X(K)$.

The geometric Bombieri--Lang conjecture is designed to take this phenomenon
into account.  A first, informal form is the following.

\begin{conjecture}[Geometric Bombieri--Lang, informal form]
Let $K/k$ be a finitely generated extension of characteristic zero, and let
$X$ be a projective variety of general type over $K$.  Then, after removing a
proper Zariski closed subset of $X$, the remaining $K$-rational points should
come from ``constant families.''
\end{conjecture}

This is only a rough form.  As in the number-field case, one can give a
more explicit statement by describing which lower-dimensional subvarieties
should be removed.  This leads to the special locus.  The precise formulation,
including the meaning of ``constant families,'' will be recalled in
Section~\ref{sec:conjectures}.

\subsection{Known results and recent progress}

The Bombieri--Lang conjecture remains very open in higher dimension.  Over
number fields, the foundational results are Faltings' theorem for curves
\cite{FaltingsMordell} and Faltings' proof of the Mordell--Lang theorem for
subvarieties of abelian varieties \cite{FaltingsLang}.  These results settle
two basic sources of examples: curves of genus at least two, and subvarieties
of abelian varieties after one removes the special locus.  Beyond these
cases, the general non-density statement for varieties of general type, and a
fortiori the stronger formulation in terms of the special locus, is still
largely open.  Finite covers of abelian varieties over number fields already
form a natural and difficult test case.  Although the expected
Bombieri--Lang statement is still open in this setting, there are weaker
arithmetic results in this direction; for instance,
Corvaja--Demeio--Javanpeykar--Lombardo--Zannier proved distribution results
for rational points on ramified covers of abelian varieties over finitely
generated fields of characteristic zero \cite{CDJZ}.

Over function fields, several important cases are known, although the general
geometric Bombieri--Lang conjecture is still far from settled.  The
function-field analogue of Mordell's conjecture for curves was proved by
Manin and Grauert \cite{Manin,Grauert}.  In higher dimension, Noguchi and
Martin-Deschamps proved finiteness results for varieties with ample cotangent
bundle \cite{NoguchiAmpleCotangent,MartinDeschamps}; later
Gillet--Roessler proved a related result over one-dimensional function fields,
in all characteristics \cite{GilletRossler}.  These theorems belong to the
hyperbolicity side of the subject: positivity of the cotangent bundle gives
strong restrictions on maps from curves, and hence on rational points over
function fields.

There are also important results for constant varieties.  In this case one
cannot hope to eliminate the constant points themselves; rather, the question
is whether all rational points are accounted for by the expected constant
contribution.  For constant targets of general type, the finiteness
theorem of Kobayashi--Ochiai for meromorphic maps is a foundational result
\cite{KobayashiOchiai}.  In the direction of Kobayashi hyperbolicity, Noguchi
proved results for hyperbolic fibre spaces over function fields
\cite{Noguchi1985}; a remaining part, related to a conjecture of Lang, was
completed in \cite{Noguchi}.

A second major line of function-field results concerns abelian varieties.  If
$X$ is a closed subvariety of an abelian variety, the geometric
Bombieri--Lang conjecture is closely related to the function-field
Mordell--Lang theorem.  Raynaud proved a hyperbolic case for subvarieties of
abelian varieties \cite{Raynaud}, and Parshin proved a related
function-field result \cite{Parshin}.  When the ambient abelian variety has
trivial $K/k$-trace, Buium and Hrushovski proved the relevant Mordell--Lang
theorem \cite{Buium,Hrushovski}.  Thus the conjectural picture is well
understood for subvarieties of abelian varieties themselves.

The recent progress discussed in this survey concerns varieties admitting
finite morphisms to abelian varieties.  Xie--Yuan first introduced partial
heights and used them to construct entire curves from sequences of rational
points, proving the algebraically hyperbolic case
\cite{XieYuanPartial}.  They then developed a more explicit construction of
linear entire curves in families of abelian varieties, which gives enough
horizontal flexibility to treat ramified covers under a trace-zero assumption
\cite{XieYuanRamified}.  The overlap of these two cases, where the cover is
hyperbolic and the ambient abelian variety has trivial trace, can also be
approached using Parshin's method \cite{Parshin}, as in work of
Bartsch--Javanpeykar \cite{BartschJavanpeykar}.  Gao subsequently removed the
trace restriction by adding a new counting argument involving the ramification
divisor \cite{GaoCounting}.

The purpose of this article is to explain the ideas behind these latest
developments, rather than only to state the final theorems.

\subsection{Vojta's dictionary as a concrete construction}

Vojta's dictionary connects Diophantine approximation with Nevanlinna theory
\cite{VojtaBook}.  One basic entry of the dictionary suggests that a sequence
of pairwise distinct rational points should be compared with an entire curve.
The main idea of Xie--Yuan is to make this comparison concrete over function
fields: from a suitable sequence of rational points, one tries to construct an
actual entire curve on a complex fiber.

Let us indicate the mechanism.  Suppose for simplicity that $K=\CC(B)$, where
$B$ is a smooth projective curve, and let $\mathcal X\to B$ be a model of
$X$.  A $K$-rational point of $X$ is a section of this family.  Under suitable
hypotheses, one may reduce to a sequence of rational points whose heights tend
to infinity.  Analytically, these heights can be interpreted as integrals, or
areas, of the corresponding sections.  If the area becomes large, then after
choosing a suitable region on the base and reparametrizing, Brody's lemma can
produce an entire curve on a fiber of $\mathcal X\to B$.

The difficulty is to control the location of such an entire curve.  There are
two kinds of bad loci.  First, there are vertical bad fibers: even if the
generic fiber has the desired general-type or hyperbolicity properties, some
special fibers of a model may degenerate.  Second, there are horizontal bad
loci.  Even inside a good fiber, there may be special subvarieties; if the
limiting entire curve lies inside them, its existence does not give the desired
contradiction.

The method discussed in this article has three main ingredients.  Partial
heights localize height growth on the base and address the vertical problem
\cite{XieYuanPartial}.  The mobility of limiting linear entire curves, coming
from the Xie--Yuan construction, gives enough horizontal flexibility in the
trace-zero ramified-cover case \cite{XieYuanRamified}.  Finally, the counting
argument in \cite{GaoCounting} treats the general trace situation by following the
ramification along the whole limiting process.  Intersections with the
ramification locus are converted into
tangencies to the branch divisor, and the proportion of such tangencies is then
controlled by Nevanlinna theory \cite{GaoCounting}.

\subsection{Organization of the article}

The rest of the article is organized around this method.  In
Section~\ref{sec:conjectures} we recall the precise conjectural framework:
the algebraic special locus, the geometric Bombieri--Lang conjecture, and the
height formulations used in the proofs.  In
Section~\ref{sec:analytic-partial} we introduce analytic special sets and
partial heights, and explain how they are used to construct entire curves from
rational points.  Section~\ref{sec:linear-entire-curves} discusses the
mobility of limiting linear entire curves.
Section~\ref{sec:gao-counting-setup} prepares the
algebraic side of the counting argument: good models, ramification, and
tangency to the branch divisor.  Section~\ref{sec:nevanlinna-expanding-discs}
records the Nevanlinna theory on expanding discs that will be used in the
counting step.  Section~\ref{sec:counting-proof} then assembles the lower and
upper estimates and completes the proof.  Finally,
Section~\ref{sec:further-directions} discusses some further directions.

\section{Statement of the geometric Bombieri--Lang conjecture}\label{sec:conjectures}

\subsection{Special loci over number fields}

The non-density form of the Bombieri--Lang conjecture, stated in
Conjecture~\ref{conj:BL-density-intro}, predicts that rational points on a
variety of general type should not be Zariski dense.  This formulation,
however, does not describe the exceptional set where rational points may
accumulate.  The special-locus formulation attempts to describe this
exceptional set geometrically.

\begin{definition}\label{def:special-locus}
Let $X$ be a projective variety over a field $K$ of characteristic zero.  The
algebraic special set $\Sp(X)$ is the Zariski closure in $X_{\Kbar}$ of the
union of the images of all non-constant rational maps
\[
  A \dashrightarrow X_{\Kbar},
\]
where $A$ ranges over abelian varieties over $\Kbar$.
\end{definition}

The intuition is that images of abelian varieties are geometric sources of
many rational points.  This leads to the following special-locus form of
the Bombieri--Lang conjecture.

\begin{conjecture}[Bombieri--Lang, special-locus form]\label{conj:BL-special}
Let $X$ be a projective variety over a number field $K$.  Then
\[
  (X\setminus\Sp(X))(K)
\]
is finite.
\end{conjecture}

This conjecture refines Conjecture~\ref{conj:BL-density-intro}: instead of
merely asserting that rational points are not Zariski dense on a variety of
general type, it predicts where the exceptional points may lie.  The
definition above is algebraic.  Later, when we discuss Vojta's dictionary and
entire curves, we will also introduce an analytic special locus defined using
holomorphic maps from $\CC$.  That analytic version is closer to Nevanlinna
theory, and the Green--Griffiths--Lang philosophy relates the analytic and
algebraic special loci.  For now, we use $\Sp(X)$ as the algebraic exceptional
locus appearing in the arithmetic conjectures.

\subsection{Constant pieces over function fields}

We now pass to function fields.  Let $K/k$ be a finitely generated extension
of characteristic zero, and assume that $k$ is algebraically closed in $K$.
The phenomenon discussed in Section~1.2 forces the function-field formulation
to allow ``constant families.''  We now make precise the birational form of
this notion.

For a projective variety $T$ over $k$, write $T_K=T\times_k K$.

\begin{definition}\label{def:birationally-constant}
Let $Z$ be a projective variety over $K$.  We say that $Z$ is
\emph{birationally constant} if there exist a projective variety $T$ over $k$ and a
birational $K$-map
\[
  \rho:T_K\dashrightarrow Z.
\]
Such a pair $(T,\rho)$ is called a \emph{constant model} of $Z$.
\end{definition}

\begin{definition}\label{def:constant-induced-points}
Let $(T,\rho)$ be a constant model of $Z$.  Let $U\subset T_K$ be the maximal
open subset on which $\rho$ is defined as a morphism
\[
  \rho^\circ:U\longrightarrow Z.
\]
The rational points of $Z$ induced by this constant model are the points in
the image of
\[
  T(k)\cap U(K)\longrightarrow U(K)\xrightarrow{\rho^\circ} Z(K),
\]
where $T(k)$ is viewed as a subset of $T_K(K)$ in the natural way.
\end{definition}

\subsection{The geometric Bombieri--Lang conjectures}

With these definitions in place, we can state the function-field versions of
Bombieri--Lang.  Throughout this subsection, $K/k$ is a finitely generated
extension of fields of characteristic zero, with $k$ algebraically closed in
$K$, and varieties are defined over $K$ unless otherwise specified.  The first
formulation is the general-type form, parallel to
Conjecture~\ref{conj:BL-density-intro}.

\begin{conjecture}[Geometric Bombieri--Lang, general-type form]\label{conj:GBL-general-type}
Let $X$ be a projective variety of general type over $K$.  If $X(K)$ is
Zariski dense in $X$, then $X$ is birationally constant.  Moreover, after
removing the rational points induced by a constant model of $X$, the remaining
$K$-rational points are not Zariski dense in $X$.
\end{conjecture}

The second formulation is more precise and parallels the special-locus form.
It says that, after removing the algebraic special locus, all remaining
rational points are accounted for by finitely many constant pieces.  The
following is the form used in \cite{XieYuanPartial}, stated with the
terminology of Definitions~\ref{def:birationally-constant} and
\ref{def:constant-induced-points}.

\begin{conjecture}[Geometric Bombieri--Lang, special-locus form]\label{conj:GBL-special}
Let $X$ be a projective variety over $K$, and let $Z$ be the Zariski closure
of $(X\setminus\Sp(X))(K)$ in $X$.  Then there are finitely many distinct
closed subvarieties $Z_1,\ldots,Z_r$ of $Z$, containing all irreducible
components of $Z$, such that each $Z_i$ is birationally constant and every
point of $(X\setminus\Sp(X))(K)$ is induced by one of these constant models.
\end{conjecture}

\subsubsection*{Examples}

The two technical features in Conjecture~\ref{conj:GBL-special} are both
necessary in general.  First, the maps from constant models should be allowed
to be birational maps rather than morphisms.  Indeed, let $K$ be a function
field of one variable over $k=\CC$, let $T$ be a smooth hyperbolic surface over
$k$, and set $Y=T_K$.  If $X$ is obtained from $Y$ by blowing up a closed point
which is transcendental over $k$, then $X$ is birationally constant, but the
birational map $T_K\dashrightarrow X$ is not a morphism.  This is
\cite[Example 2.11]{XieYuanPartial}.

Second, the finite family $\{Z_1,\ldots,Z_r\}$ may have to contain more
subvarieties than just the irreducible components of $Z$.  Let $C_0$ be a
smooth projective curve of genus $>1$ over $k=\CC$, let $K=k(C_0)$, and set
$C=C_0\times_k K$.  The set
\[
  \Sigma=C(K)\setminus C_0(k)
\]
is finite and nonempty.  Now take $X=C\times_K C$.  Then $X$ is constant and
algebraically hyperbolic, and
\[
  X(K)=C(K)\times C(K)
\]
is the union of
\[
  C_0(k)\times C_0(k),\qquad
  C_0(k)\times \Sigma,\qquad
  \Sigma\times C_0(k),\qquad
  \Sigma\times\Sigma .
\]
Thus the irreducible component $X$ itself does not account for all rational
points via $X_0(k)$.  One also needs the constant curves
$P\times_K C$ and $C\times_K P$ for $P\in\Sigma$.  This explains why the
finite set $\{Z_1,\ldots,Z_r\}$ in the conjecture may contain lower-dimensional
constant pieces in addition to the irreducible components of $Z$; see
\cite[Example 2.12]{XieYuanPartial}.

\subsection{Heights over function fields}

For the proofs discussed in this article, it is useful to translate the
conjecture into height language.  We recall only the function-field picture
needed later.  Suppose that $K=k(B)$, where $B$ is a smooth projective curve
over $k$.  Let $\mathcal X\to B$ be a projective model of $X$, and let
$\mathcal L$ be a relatively ample line bundle on $\mathcal X$.  A point
$x\in X(K)$ is represented by a section $\sigma_x:B\to\mathcal X$, and its
height is, up to the usual bounded ambiguity,
\[
  h_{\mathcal L}(x)=\deg \sigma_x^*\mathcal L .
\]
When $k=\CC$, this degree can also be viewed analytically as an area of the
section.

Over number fields, bounded height is closely tied to finiteness by
Northcott's theorem.  Over function fields of characteristic zero, this is no
longer true: constant points may form infinite families of bounded height.
This is one reason the height formulation is natural in the geometric setting.

\begin{conjecture}[Geometric Bombieri--Lang, height form]\label{conj:GBL-height}
Let $X$ be a projective variety of general type over $K$.  Then there exists a
nonempty open subset $U\subset X$ such that the height of the points in
$U(K)$ is bounded.
\end{conjecture}

Equivalently, one can use \emph{generic sequences}.  A sequence
$\{x_n\}_{n\ge 1}\subset X(K)$ is called generic if every infinite subsequence
is Zariski dense in $X$.

\begin{conjecture}[Generic-sequence form]\label{conj:GBL-generic}
Let $X$ be a projective variety of general type over $K$.  Then every generic
sequence of $K$-rational points in $X$ has bounded height.
\end{conjecture}

This is the form most directly connected with the analytic constructions
below.  To prove the conjecture by contradiction, one starts with a generic
sequence whose heights tend to infinity and tries to construct entire curves,
or to perform a counting argument, from that sequence.

\subsection{Relations among the formulations}

The different formulations above are not independent.  The
special-locus form is a more explicit formulation.  Under the expected
relation between general type and the special locus, such as the
Green--Griffiths--Lang picture in the analytic setting, this explicit form
implies the general-type form: once $\Sp(X)$ is a proper closed subset,
controlling rational points outside $\Sp(X)$ gives the desired non-density
statement.  The height form and the generic-sequence form are also equivalent,
and they are often the most convenient versions for proofs.

For the class of varieties admitting finite morphisms to abelian varieties,
the needed relation between general type and the special locus is available.
Consequently, proving the height or generic-sequence version in this class
implies the corresponding geometric Bombieri--Lang statements.  This is the
framework in which the results of Xie--Yuan and Gao will be discussed.

Finally, the algebraic special locus is not the only special set that will
appear.  When entire curves enter the argument, an analytic special locus
defined using holomorphic maps from $\CC$ becomes natural.  We postpone its
definition to the next part of the article, where the analytic side of the
argument enters.

\subsection{Varieties finite over abelian varieties}

We now state the main result discussed in this survey.  It will be useful to
record it in two forms.  The first is the geometric form, expressed in terms of
Conjecture~\ref{conj:GBL-special}.

\begin{theorem}[Xie--Yuan; Gao, geometric form]\label{thm:finite-cover-geometric}
Let $X$ be a projective variety over $K$ admitting a finite morphism to an
abelian variety over $K$.  Then Conjecture~\ref{conj:GBL-special} holds for
$X$.
\end{theorem}

Equivalently, in the general-type case one can use the generic-sequence
formulation.  For the proof strategy discussed in the rest of the article, we
use the following curve-function-field form.

\begin{theorem}[Xie--Yuan; Gao, generic-sequence form]\label{thm:finite-cover-generic}
Let $K=\CC(B)$ be the function field of a smooth projective complex curve $B$,
arising as the auxiliary curve in the standard reduction.  Let $X$ be a
projective variety of general type over $K$ admitting a finite morphism to an
abelian variety over $K$.  Then every generic sequence in $X(K)$ has bounded
height.
\end{theorem}

The equivalence of these two forms in the finite-cover setting follows from
the comparison above and the known description of the special locus in this
class.  The hyperbolic finite-cover case was proved by Xie--Yuan using partial
heights and entire curves \cite{XieYuanPartial}.  The ramified-cover case with
trivial trace was then treated by Xie--Yuan using linear entire curves
\cite{XieYuanRamified}.  Gao completed the general finite-cover case by a new
counting argument \cite{GaoCounting}.  In the rest of the article we will
mainly explain the proof strategy for the generic-sequence form.

Thus, whenever the large field $K$ appears in the proof-oriented sections
below, it should be understood as this auxiliary curve-function field
$K=\CC(B)$ unless stated otherwise.

\section{Analytic special sets and partial heights}\label{sec:analytic-partial}

We begin with the analytic side of the construction.  The main technical point
of this section is the use of partial heights to localize the Brody
construction on the base curve.  Throughout this section we work in the
curve-function-field setting $K=\CC(B)$ introduced in the generic-sequence
formulation of Theorem~\ref{thm:finite-cover-generic}.

\subsection{The analytic special locus}

The algebraic special locus in Definition~\ref{def:special-locus} is the
exceptional set appearing in the arithmetic conjectures.  There is also a
more immediate analytic analogue.  It is motivated by the basic entry of
Vojta's dictionary which compares an infinite sequence of rational points with
an entire curve.

\begin{definition}\label{def:analytic-special-locus}
Let $X$ be a projective variety over $\CC$.  The \emph{analytic special set}
$\Span(X)$ is the Zariski closure in $X$ of the union of the images of all
non-constant holomorphic maps
\[
  f:\CC\longrightarrow X(\CC).
\]
We say that $X$ is \emph{Brody hyperbolic} if $\Span(X)=\emptyset$.
\end{definition}

This definition records exactly where entire curves can lie.  It is therefore
the analytic counterpart of the algebraic special locus.

\subsection{Lang's conjecture and Green--Griffiths--Lang}

Lang conjectured that the algebraic and analytic special sets should agree:
\[
  \Span(X)=\Sp(X)
\]
for projective varieties over $\CC$.  The Green--Griffiths--Lang conjecture is
a closely related statement: it predicts that a projective variety is of
general type precisely when its analytic special set is a proper closed subset.
Both conjectures are far beyond current methods in general.  In the class
needed in this survey, however, the required comparison is known.

\begin{theorem}[Ueno--Kawamata--Noguchi--Winkelmann--Yamanoi]\label{thm:special-comparison}
Let $X$ be a projective variety over $\CC$ admitting a finite morphism to an
abelian variety.  Then
\[
  \Span(X)=\Sp(X).
\]
Moreover, $X$ is of general type if and only if $\Sp(X)\ne X$.
\end{theorem}

We use this theorem only in the form stated above; see \cite{Yamanoi} and
\cite[Theorem 4.1 and Corollary 4.3]{XieYuanPartial}.  It is what allows one
to translate a holomorphic curve produced by an analytic argument into
information about the algebraic special locus.

The relevance to the geometric Bombieri--Lang conjecture is the following.
If a sequence of rational points over a function field produces an entire
curve on a complex fiber, then the image of that curve lies in $\Span(X)$.  In
a hyperbolic situation, where $\Span(X)=\emptyset$, this gives an immediate
contradiction.  In general, one needs more information about the location of
the entire curve.  This is the source of the horizontal difficulty discussed
in the introduction.

\subsection{From sections to areas}

We now return to a function field $K=\CC(B)$, where $B$ is a smooth projective
complex curve.  Let $X$ be a projective variety over $K$, and choose a
projective model
\[
  \pi:\mathcal X\longrightarrow B
\]
with generic fiber $X$.  A point $x\in X(K)$ is a section
$\sigma_x:B\to\mathcal X$.  After choosing a relatively ample line bundle
$\mathcal L$ on $\mathcal X$, the height of $x$ is the degree
\[
  h_{\mathcal L}(x)=\deg \sigma_x^*\mathcal L .
\]
Analytically, after choosing a smooth metric, this degree is represented by an
area integral
\[
  \int_B \sigma_x^* c_1(\mathcal L).
\]

Thus a sequence of rational points whose heights tend to infinity becomes a
sequence of sections whose areas tend to infinity.  This is the starting point
for turning rational points into entire curves.  If one can find a place on
the base where a definite amount of area concentrates, then one can rescale
the section near that place.  Brody's lemma then gives a subsequential limit,
which is a non-constant holomorphic map from $\CC$ to a complex fiber of
$\pi$.

\subsection{Partial heights on projective varieties}

The total height alone does not say where the area is located on the base.
This is the reason for introducing \emph{partial heights}.  Let
$\pi:\mathcal X\to B$ be a projective model of $X$ and let
$D\subset B$ be a measurable subset whose closure is compact.  Let $\omega$ be
a real $(1,1)$-form on a neighborhood of $\pi^{-1}(\overline D)$.  For a
section $x:B\to\mathcal X$, the partial height attached to the pair $(D,\omega)$
is
\[
  h_{(D,\omega)}(x)=\int_D x^*\omega .
\]
More generally, for a closed point of $X$ one integrates over the corresponding
multi-section and divides by the degree.  If $D=B$ and $\omega$ represents the
first Chern class of a line bundle $\mathcal L$ extending an ample line bundle
on $X$, this recovers the usual full height.

Partial heights are bounded above by full heights.  More precisely, if
$h_L$ is a Weil height attached to an ample line bundle $L$ on $X$, then for
any partial height $h_{(D,\omega)}$ there are constants $c_1,c_2>0$ such that
\[
  |h_{(D,\omega)}(x)|\le c_1 h_L(x)+c_2
\]
for all $x\in X(K)$.  This is the easy direction.  The difficult point is
whether growth of the full height forces growth of a given positive partial
height.

\begin{conjecture}[Xie--Yuan non-degeneracy conjecture]\label{conj:nondeg}
Let $K=\CC(B)$, where $B$ is a smooth projective curve.  Let $X$ be a
projective variety over $K$ containing no rational curve.  Let $h_L$ be a Weil
height attached to an ample line bundle on $X$, and let $h_{(D,\omega)}$ be a
partial height attached to a strictly positive pair $(D,\omega)$, where $D$ is
an open disc in $B$.  If $\{x_n\}\subset X(K)$ satisfies
\[
  h_L(x_n)\longrightarrow \infty,
\]
then
\[
  h_{(D,\omega)}(x_n)\longrightarrow \infty .
\]
\end{conjecture}

This is \cite[Conjecture 2.3]{XieYuanPartial}.  It says that positive partial
heights are not merely dominated by full heights; for varieties without
rational curves, they should detect the same unbounded sequences.
The assumption excluding rational curves is necessary.  Xie--Yuan give an
example for $X=\mathbb P^1_K$: one can choose sections whose usual heights go
to infinity, while after a suitable rescaling their partial heights over a
fixed disc remain bounded \cite[Example 2.4]{XieYuanPartial}.

\subsection{Canonical partial heights on abelian varieties}

The non-degeneracy conjecture is proved in \cite{XieYuanPartial} for varieties
finite over abelian varieties.  The key input is a canonical version of partial
height on abelian varieties.

Let $A$ be an abelian variety over $K=\CC(B)$, and let $L$ be a symmetric ample
line bundle on $A$.  After shrinking $B$, assume that $A$ extends to an
abelian scheme $\mathcal A\to B$.  Let $\omega$ be the Betti form on
$\mathcal A$ associated to $L$.  For a section $s\in A(K)=\mathcal A(B)$ and
a measurable subset $D\subset B$, define the \emph{canonical partial height}
by
\[
  \widehat h_{(D,\omega)}(s)=\int_D s^*\omega .
\]
If $D=B$, then a theorem of Gauthier--Vigny \cite{GauthierVigny} identifies
this integral with the Neron--Tate height:
\[
  \widehat h_L(s)=\int_B s^*\omega .
\]
In particular,
\[
  0\le \widehat h_{(D,\omega)}(s)\le \widehat h_L(s).
\]
Like the Neron--Tate height, the canonical partial height is quadratic:
\[
  \widehat h_{(D,\omega)}([m]s)=m^2\widehat h_{(D,\omega)}(s)
\]
and it satisfies the parallelogram identity.  Thus it induces a positive
semidefinite quadratic form on $A(K)\otimes_{\mathbb Z}\mathbb R$.

The crucial result is that, when $D$ is a disc, this quadratic form is
non-degenerate modulo the constant part.

\begin{theorem}[Xie--Yuan]\label{thm:canonical-partial-nondeg}
Let $D\subset B$ be an open disc, and write $A^{K/\CC}$ for the $K/\CC$-trace
of $A$.  Then the canonical partial height $\widehat h_{(D,\omega)}$ is
invariant under translation by $A^{K/\CC}(\CC)$, and it induces a positive
definite quadratic form on
\[
  \bigl(A(K)/A^{K/\CC}(\CC)\bigr)\otimes_{\mathbb Z}\mathbb R.
\]
Consequently, there is a constant $\epsilon>0$ such that
\[
  \widehat h_{(D,\omega)}(s)\ge \epsilon\,\widehat h_L(s)
\]
for all $s\in A(K)$.
\end{theorem}

This is \cite[Theorem 3.5]{XieYuanPartial}.  Together with the comparison of
ordinary partial heights under finite morphisms, it gives the following
consequence, which is the form used later.

\begin{corollary}[Xie--Yuan]\label{cor:nondeg-finite-abelian}
Let $K=\CC(B)$, where $B$ is a smooth projective curve.  Let $X$ be a
projective variety over $K$ admitting a finite morphism to an abelian variety.
Then the non-degeneracy conjecture, Conjecture~\ref{conj:nondeg}, holds for
$X$.
\end{corollary}

This is \cite[Theorem 3.6]{XieYuanPartial}.

\subsection{The hyperbolic finite-cover case}

We now explain how the preceding ingredients prove the hyperbolic case of the
finite-cover theorem.  Let $X$ be a projective variety over $K=\CC(B)$
admitting a finite morphism to an abelian variety.  Assume that $X$ is
algebraically hyperbolic in the sense that $\Sp(X)=\emptyset$.\footnote{Following
\cite{XieYuanPartial}, we use the term algebraically hyperbolic for the
condition $\Sp(X)=\emptyset$.  This condition is also called \emph{groupless}
in parts of the literature; it should not be confused with Demailly's notion
of algebraic hyperbolicity.}

Choose a projective model $\pi:\mathcal X\to B$ of $X$, after replacing $B$ by
a nonempty open subset if necessary.  Since $X$ is algebraically hyperbolic,
the very general fibers of $\pi$ are also algebraically hyperbolic.  More
precisely, after removing a countable union of proper closed subsets from the
base, the specialization of the algebraic special set remains empty on the
fiber.  For such a very general point $b\in B(\CC)$, the fiber $\mathcal X_b$
is a complex projective variety finite over an abelian variety and satisfies
\[
  \Sp(\mathcal X_b)=\emptyset .
\]
By Theorem~\ref{thm:special-comparison}, the algebraic and analytic special
sets agree for varieties finite over abelian varieties.  Hence
\[
  \Span(\mathcal X_b)=\emptyset ,
\]
so $\mathcal X_b$ is Brody hyperbolic.

We use the standard openness of Brody hyperbolicity in proper holomorphic
families of compact complex analytic spaces.  Thus, after choosing $b$ as
above and shrinking around it, we may find a small disc $D\subset B$ such that
all fibers $\mathcal X_t$, $t\in D$, are Brody hyperbolic.

Suppose now, toward a contradiction, that the generic-sequence form,
Conjecture~\ref{conj:GBL-generic}, fails.
Then there is a generic sequence $\{x_n\}\subset X(K)$ whose heights tend to
infinity.  By Corollary~\ref{cor:nondeg-finite-abelian}, the non-degeneracy
conjecture holds for $X$; hence the partial height of $x_n$ over the disc $D$
also tends to infinity.  In analytic terms, the corresponding sections
$\sigma_{x_n}:B\to\mathcal X$ have unbounded area over $D$.

Applying Brody's lemma to the restrictions $\sigma_{x_n}|_D$, one obtains,
after reparametrization and passing to a subsequence, a non-constant entire
curve
\[
  \CC\longrightarrow \mathcal X_t
\]
for some $t\in D$.  But every fiber over $D$ is Brody hyperbolic, so this is
impossible.  This contradiction proves the hyperbolic finite-cover case.

This argument also indicates the limitation of partial heights alone.  Partial
heights solve the vertical problem: they allow the Brody construction to be
localized over a disc where the fibers are hyperbolic.  In the non-hyperbolic
case, however, there may be horizontal bad loci, namely special subvarieties
inside the good fibers.  A Brody limit may lie entirely inside such a locus.
This leads to the need for a mobility statement for limiting linear entire
curves.

\section{Mobility of limiting linear entire curves}\label{sec:linear-entire-curves}

Partial heights localize the Brody construction, but the resulting entire
curve may still lie in a horizontal bad locus.  The next ingredient is the
mobility of limiting linear entire curves in fibers of abelian schemes.  The
point is not merely that the limiting curves are linear; rather, their
limiting positions have enough freedom to avoid horizontal bad loci.  We
recall the part of the Xie--Yuan construction needed later
\cite{XieYuanRamified}.  We
keep the same curve-function-field setting $K=\CC(B)$ as in
Theorem~\ref{thm:finite-cover-generic}.

\subsection{The trace of an abelian variety}

We first recall the trace, since it measures the constant part of an abelian
variety over a function field.  Let $K/k$ be a field extension as above, and
let $A$ be an abelian variety over $K$.  The \emph{$K/k$-trace} of $A$ is an
abelian variety $A^{K/k}$ over $k$, together with a $K$-homomorphism
\[
  \operatorname{tr}:A^{K/k}\times_k K\longrightarrow A,
\]
which is universal among homomorphisms from constant abelian varieties to
$A$.  In characteristic zero, the trace map is a closed immersion.  Thus one
may regard $A^{K/k}$ as the constant part of $A$.

The case $A^{K/k}=0$ is often called the \emph{trace-zero} case.  In that
case the sections of an abelian scheme have no constant directions in which
they can drift.  A constant factor may restrict the freedom of the limiting
linear entire curves.  The counting argument supplies the additional input
needed to compensate for this possible loss of freedom.

\subsection{Linear entire curves in abelian varieties}

Let $A$ be a complex abelian variety.  Every vector $v\in \operatorname{Lie}(A)$
defines a holomorphic map
\[
  \CC\longrightarrow A,\qquad z\longmapsto x+\exp_A(zv),
\]
where $x\in A$.  We call such a map a \emph{linear entire curve}.  Its image is
a translate of a one-parameter complex subgroup, and its Zariski closure is a
translate of an abelian subvariety of $A$.

In families, this simple construction becomes a way to extract controlled
entire curves from sections.  Let $\pi:\mathcal A\to B$ be an abelian scheme
over a smooth complex curve, with generic fiber $A/K$ for $K=\CC(B)$.  A
section $s\in A(K)=\mathcal A(B)$ may be studied by looking at its behavior in
Betti coordinates.  The derivative of the Betti coordinates gives a transfer
map
\[
  \delta(s,\cdot):T_bB\longrightarrow \operatorname{Lie}(\mathcal A_b).
\]
This map records the infinitesimal direction in which the section moves inside
the fiber when one rescales the base near $b$.

The following non-degeneracy property is one of the inputs from
\cite{XieYuanRamified}.  If $s$ is nonzero modulo the constant part, then
$\delta(s,\cdot)$ is nonzero for all $b$ outside a set of measure zero.  Thus,
for very general $b$, a nonconstant sequence of sections gives a genuine
linear direction in the fiber.

\subsection{Limits of reparametrized sections}

The basic construction is the following.  Let $\{s_n\}\subset A(K)$ be a
sequence of sections.  After rescaling by positive real numbers $\ell_n\to
\infty$, suppose that the images of $\ell_n^{-1}s_n$ in the real vector space
obtained from the Mordell--Weil group converge to a nonzero vector
$s_\infty$.  Let $b_n\to b$ and suppose that $s_n(b_n)$ converges to a point
$x\in \mathcal A_b$.  If the transfer map
\[
  \delta(s_\infty,\cdot):T_bB\longrightarrow \operatorname{Lie}(\mathcal A_b)
\]
is nonzero, then a suitable reparametrization of the sections $s_n$ converges,
uniformly on compact subsets of $\CC$, to the linear entire curve in
$\mathcal A_b$ passing through $x$ in the direction
$\delta(s_\infty,v_b)$.

This is the content of the Xie--Yuan convergence theorem
\cite[Theorem 3.1]{XieYuanRamified}; see also the formulation recalled in
\cite[Theorem 2.3]{GaoCounting}.  The important point for the survey is that
the limiting entire curve is not arbitrary: its direction and its Zariski
closure can be described in terms of $s_\infty$ and the transfer map.

\subsection{The mobility statement}

The construction above also gives many possible limit points in a very general
fiber.  We will use the following consequence of the Xie--Yuan construction,
in the quotient-height form recorded in \cite[Corollary 2.5]{GaoCounting}.

\begin{proposition}[Xie--Yuan, in the form used by Gao]\label{prop:many-linear-limits}
Let $\pi:\mathcal A\to B$ be an abelian scheme over a smooth complex curve
with generic fiber $A/K$, and let $\{s_n\}\subset A(K)$ be a sequence of
sections.  Assume that the image of $\{s_n\}$ in every nontrivial quotient
abelian variety of $A$ has height tending to infinity.  Assume also the
standard endomorphism condition appearing in the Xie--Yuan construction, namely
that all geometric endomorphisms of $A$ are defined over $K$.  Then, outside a
subset of $B(\CC)$ of measure zero, for every $b$ there is a Zariski dense
subset $U_b\subset \mathcal A_b$ such that every point of $U_b$ occurs as a
limit point of the images of the sections $s_n(B)$.
Moreover, these limit points arise from reparametrizations whose limits are
linear entire curves in the fiber $\mathcal A_b$.
\end{proposition}

The formulation above is intentionally qualitative.  The proof combines the
description of Zariski closures of linear leaves, the non-degeneracy of the
transfer map, and the convergence theorem for reparametrized sections.  This
is the mobility property that will later allow the construction to avoid
horizontal bad loci.  The quotient-height hypothesis is the essential
non-degeneracy condition.  In the trace-zero situation, for a generic sequence,
Northcott's property on the Mordell--Weil group modulo the trace implies the
corresponding unboundedness after passing to the quotients which occur in the
standard reductions.  This hypothesis is used in the proof by induction on the
number of simple factors of $A$: when $A$ is simple there is no proper quotient
to control, while in general one takes the smallest abelian subvariety
containing the limiting direction and applies the induction hypothesis to the
quotient.

The condition is also necessary.  For instance, if
$A=A_0\times A_1$ and the $A_0$-coordinate of $s_n$ is fixed while the
$A_1$-coordinate has height tending to infinity, then the limiting positions
may be fixed in the $A_0$-direction and need not be Zariski dense in the whole
fiber.  In the application to the reduction in
Proposition~\ref{prop:gao-generic-reduction}, the quotient-height condition
and the endomorphism condition are achieved after the standard quotient
reductions and base change.

Let us also indicate how this mobility is used in the trace-zero
argument of Xie--Yuan.  The limiting linear entire curve is first visible in a
fiber of the abelian scheme $\mathcal A$.  Since the original sections come
from $X$, one applies Brody's lemma to the lifted sections on a model of $X$;
after passing to a subsequence, the limiting curve in the abelian fiber lifts
to an entire curve in the corresponding fiber of $X$.  The mobility statement
then allows one to choose the limiting position so that the lifted entire curve
is not contained in the special locus, giving the contradiction in the
trace-zero case.  This is different from the counting argument below, where
the final contradiction no longer comes from avoiding $\Sp(X)$, but from
ramification and Nevanlinna estimates.

\section{Ramification and good models}\label{sec:gao-counting-setup}

The preceding sections explain how rational points of large height give rise
to limiting linear entire curves with enough mobility.  In the general trace
case this mobility alone is not enough, and the new input is to use the
ramification of the finite morphism.

The counting argument starts from a simple geometric observation.  The
points under consideration are sections of $X$, not arbitrary sections of the
abelian variety.  When such a section meets the ramification divisor of the
finite morphism, its image on the abelian side meets the branch locus
non-transversally; in other words, the projected section is tangent to the
reduced branch divisor.  Thus ramification converts intersections on $X$ into
tangencies on the abelian side.

The proof then compares two opposite estimates for these tangencies.  On the
algebraic side, the positivity of the ramification divisor forces tangencies
to occur with positive proportion.  This positivity ultimately comes from the
canonical divisor of $X$: after passing to a suitable model, the ramification
divisor carries the general-type positivity of $X$.  On the analytic side,
Nevanlinna theory applied to the limiting linear entire curves shows that such
tangencies should have asymptotic proportion zero.  The contradiction between
these two estimates is the heart of the counting argument.

In this section we discuss the algebraic-geometric preparation needed for the
first estimate.  We reduce to a finite surjective morphism, arrange a good
model on which ramification can be measured, and record the tangency
mechanism.  The actual counting functions and the Nevanlinna upper bound will
be treated separately in the next section.

\subsection{Reduction to the finite surjective case}

We first reduce Theorem~\ref{thm:finite-cover-generic} to the following core
boundedness statement.

\begin{proposition}[Key boundedness statement]\label{prop:gao-generic-reduction}
Let $K=\CC(B)$, where $B$ is a smooth projective complex curve.  Let $X$ be a
projective variety of general type over $K$, let $A$ be an abelian variety over
$K$, and let $T$ be a projective variety over $\CC$.  Suppose that there is a
finite surjective morphism
\[
  f:X\longrightarrow A\times_K T_K .
\]
Let $\{x_n\}\subset X(K)$ be a generic sequence.  Write
\[
  f(x_n)=(s_n,t_n)\in A(K)\times T_K(K).
\]
If each $t_n$ is induced by a point of $T(\CC)$, then the heights of the
points $x_n$ are bounded.
\end{proposition}

\subsubsection*{Sketch of the reduction}

Assume that Theorem~\ref{thm:finite-cover-generic} fails.  Then there is a
generic sequence of rational points whose heights tend to infinity.  Applying
the trace and quotient reductions to the abelian factor, one separates the
genuinely varying abelian part from the constant part.  The latter is recorded
by an auxiliary projective variety $T$ over $\CC$.  In this way one obtains a
finite morphism
\[
  f:X\longrightarrow A\times_K T_K
\]
such that the $T$-coordinates of the sequence come from $T(\CC)$.  Thus the
constant $T$-coordinate records the constant part produced by the reduction;
boundedness of the whole sequence is still the assertion to be proved.

It remains to explain why one may pass to the surjective case.  Let
$Z\subset A\times_KT_K$ be the image of $X$.  After replacing $T$ by the
Zariski closure of the points $t_n$, we may assume that the sequence $\{t_n\}$
is generic in $T$.  Since the sections $(s_n,t_n)$ lie in $Z$, the projection
$Z\to T_K$ is then dominant.  Suppose that $Z$ is still a proper closed
subvariety of $A\times_KT_K$.

Choose models over $B$ and let
\[
  \mathcal Z\subset \mathcal A\times_\CC T
\]
be the closure of $Z$.  By Raynaud--Gruson's flattening theorem
\cite{RaynaudGruson}, after replacing $T$ by a birational modification
$T'\to T$, the strict transform $\mathcal Z'\subset \mathcal A\times_\CC T'$
may be assumed flat over $T'$.  The sections $(s_n,t_n)$ lift, after passing
to a subsequence, to sections $(s_n,t'_n)$ contained in $\mathcal Z'$, with
$t'_n\to t'_\infty$ in $T'(\CC)$.

Because $Z$ is proper in $A\times_KT_K$ and $\mathcal Z'\to T'$ is flat, the
fiber $\mathcal Z'_{b,t'_\infty}$ is a proper closed subset of
$\mathcal A_b\times\{t'_\infty\}$ for very general $b\in B(\CC)$.  On the
other hand, the mobility statement for limiting linear entire curves gives a
Zariski dense set of possible limit points of the sections $s_n(B)$ in
$\mathcal A_b$.  We may therefore choose such a limit point $y\in\mathcal A_b$
with
\[
  (y,t'_\infty)\notin \mathcal Z'_{b,t'_\infty}.
\]
This is impossible, because all sections $(s_n,t'_n)$ are contained in the
closed subset $\mathcal Z'$, so their limits must lie in the limiting fiber
$\mathcal Z'_{b,t'_\infty}$.  Thus the non-surjective case is eliminated by
the mobility of limiting linear entire curves.

We are therefore reduced to the finite surjective situation of
Proposition~\ref{prop:gao-generic-reduction}.  The proposition then gives
bounded height, contradicting the choice of the original sequence.  Hence
Proposition~\ref{prop:gao-generic-reduction} implies
Theorem~\ref{thm:finite-cover-generic}.  \hfill$\square$

We now turn to the proof of Proposition~\ref{prop:gao-generic-reduction}.
After replacing $X$ by its normalization and passing to a subsequence if
necessary, we may assume that $X$ is normal.  Suppose, for contradiction, that
the heights of $x_n$ tend to infinity.  After further quotient reductions on
the abelian factor, we may assume that the image of the sequence in every
non-trivial quotient of the varying abelian part has height tending to
infinity.  This is the condition needed to apply the mobility statement for
limiting linear entire curves.  Hence, on very general fibers, the sections $s_n$
have many possible limiting positions.  The remaining algebraic work is to
place the ramification of $f$ into a usable form.

\subsection{Good models and ramification}

We keep the notation of Proposition~\ref{prop:gao-generic-reduction}.  Let
$R\subset X$ be the ramification divisor of the finite morphism
$f:X\to A\times_KT_K$.  Since the target is smooth and $X$ is normal, this
divisor is first defined on the smooth locus of $X$ and then extended by
closure as a Weil divisor on $X$.  Let $D=f(R)$ be the corresponding branch
divisor on $A\times_KT_K$.

The first reduction concerns the position of ramification over the constant
factor $T$.  A component of the branch divisor could contain a whole fiber
$A\times\{t\}$.  Such a component is not useful for the later tangency count
on the fiber over $t$, because it does not define a genuine divisor on that
fiber.  Gao removes this difficulty by using weak semistable reduction of
Abramovich--Karu \cite{AbramovichKaru}.  After replacing $T$ by a generically
finite modification $T'\to T$ and taking the corresponding strict transform
and normalization of $X$, one may arrange that the relevant branch divisor
does not contain any fiber $A\times\{t'\}$.  Equivalently, the relevant
ramification is horizontal over the constant factor.

After this reduction, one still has to measure the positivity of ramification.
Although $X$ is normal, the ramification divisor on a naive model need not be
$\mathbb Q$-Cartier.  One therefore passes to the canonical model of
$X$, in the sense of the minimal model program and the theorem of
Birkar--Cascini--Hacon--McKernan \cite{BCHM}.  On
this model the canonical class is ample and $\mathbb Q$-Cartier.  The original
finite morphism $X\to A\times_KT_K$ is still needed for the branch divisor and
the tangency statement.  The canonical model has a different role: it is the
place where the positivity comparison is made.  Since the target
$A\times_KT_K$ has trivial canonical bundle along the abelian fibers, the
ramification formula on the canonical model
\[
  K_{X_{\mathrm{can}}}
  \sim_{\mathbb Q}
  (f_{\mathrm{can}})^*K_{A\times_KT_K}
  + R_{\mathrm{can}}
\]
says that the positivity of $K_{X_{\mathrm{can}}}$ is carried by the
ramification divisor.  The ampleness of $K_{X_{\mathrm{can}}}$ is essential
here; a merely big divisor would leave base-locus issues that are not suitable
for the later uniform counting estimates.  Thus, after the reductions,
ramification is not a small error term: it is the source of the general-type
positivity.  This is the algebraic input that will later become a lower bound
in the counting argument.

For the integral models used below, it is convenient to resolve the
canonical model once and for all.  Thus the positivity is measured on
$X_{\mathrm{can}}$, while the actual divisors used for intersections and
tangencies are pulled back to a smooth birational model
$\widetilde X\to X_{\mathrm{can}}$.

\begin{proposition}[Good model package]\label{prop:good-model-package}
After passing to a subsequence and replacing $(X,T,f)$ by the modifications
above, we may assume that the following properties hold.
\begin{enumerate}
\item The variety $X$ is normal and there is a finite surjective morphism
\[
  f:X\longrightarrow A\times_KT_K .
\]
If $D\subset A\times_KT_K$ is the branch divisor of $f$, then $D$ contains no
whole fiber $A\times\{t\}$ over a point $t\in T(\CC)$.  The constant
coordinates of the sequence still come from points of $T(\CC)$.
\item There is a canonical model $X_{\mathrm{can}}$ in the sense of
BCHM and a birational morphism
\[
  \varpi:X_{\mathrm{can}}\longrightarrow X,
\]
so that the composite
\[
  f_{\mathrm{can}}=f\circ\varpi:
  X_{\mathrm{can}}\longrightarrow A\times_KT_K
\]
is generically finite.  The divisor $K_{X_{\mathrm{can}}}$ is ample and
$\mathbb Q$-Cartier, and the ramification divisor $R_{\mathrm{can}}$ of
$f_{\mathrm{can}}$ is $\mathbb Q$-Cartier.  Moreover
\[
  K_{X_{\mathrm{can}}}
  \sim_{\mathbb Q}
  f_{\mathrm{can}}^*K_{A\times_KT_K}+R_{\mathrm{can}} .
\]
Since $K_{A\times_KT_K}$ is trivial along the abelian fibers,
$R_{\mathrm{can}}$ is relatively ample over $T_K$.
\item Choose a resolution
\[
  \mu:\widetilde X\longrightarrow X_{\mathrm{can}}
\]
and set $\widetilde f=f_{\mathrm{can}}\circ\mu$.  Let
$\widetilde R=\mu^*R_{\mathrm{can}}$.  After multiplying by a fixed positive
integer if necessary, the divisors used on $\widetilde X$ are Cartier.
\item If $R$ denotes the ramification divisor of the finite morphism $f$, then
\[
  R=\varpi(R_{\mathrm{can}})
\]
as closed subsets of $X$.  This is \cite[Lemma 3.13]{GaoCounting}; it follows
from the negativity lemma, since no extra exceptional component of
$R_{\mathrm{can}}$ can lie over $X\setminus R$.
\item After choosing projective models over $B$, there is a smooth projective
model $\widetilde{\mathcal X}$ of $\widetilde X$ and an extension
\[
  \widetilde{\mathcal X}\longrightarrow \mathcal A\times_\CC T
\]
of $\widetilde f$.  Let
$\widetilde{\mathcal R}\subset\widetilde{\mathcal X}$ be the Zariski
closure of $\widetilde R$, and let
$\mathcal D\subset\mathcal A\times_\CC T$ be the Zariski closure of its branch
divisor.  These are the divisors used in the tangency and counting arguments.
\end{enumerate}
In diagrammatic form, the generic fibers fit into
\[
\begin{array}{ccccc}
\widetilde X & \xrightarrow{\ \mu\ } & X_{\mathrm{can}}
& \xrightarrow{\ \varpi\ } & X \\[3mm]
{\scriptstyle \widetilde f}\searrow && {\scriptstyle f_{\mathrm{can}}}\downarrow
&& \swarrow{\scriptstyle f} \\[1mm]
& & A\times_KT_K & &
\end{array}
\]
\end{proposition}

In the sequel, all references to $X,T,f$ are made after this replacement.  A
rational point in the resulting sequence gives, after lifting to the chosen
model when necessary, a section
\[
  \widetilde x:B\longrightarrow\widetilde{\mathcal X},
\]
and its image gives a section $(s,t):B\to\mathcal A\times_\CC T$.

\subsection{Tangency from ramification}

We now work in the good model package of
Proposition~\ref{prop:good-model-package}.  The elementary geometric
observation behind the counting argument is that ramification forces
non-transversality on the abelian side.

\begin{proposition}[Tangency along the branch divisor]\label{prop:tangency-branch}
Let $\widetilde x:B\to\widetilde{\mathcal X}$ be a section and write its image in
$\mathcal A\times_\CC T$ as $(s,t)$.  Suppose that, for some $b\in B(\CC)$,
the section $\widetilde x$ meets the ramification divisor
$\widetilde{\mathcal R}$ at $b$.
Then the section
\[
  (s,t):B\longrightarrow \mathcal A\times_\CC T
\]
is tangent at $(s(b),t(b))$ to the reduced branch divisor
$\mathcal D_{\mathrm{red}}$.
Equivalently, the image of
\[
  T_bB\longrightarrow T_{(s(b),t(b))}(\mathcal A\times_\CC T)
\]
is contained in $T_{(s(b),t(b))}\mathcal D_{\mathrm{red}}$.
\end{proposition}

This is the geometric content of \cite[Proposition 1.7]{GaoCounting};
the form used here follows from \cite[Proposition 3.19 and
Corollary 3.20]{GaoCounting}.  Intuitively, ramification means that the
differential of the finite map drops rank.  At a ramification point, the image
of the tangent direction of the lifted section therefore lies in the tangent
space of the branch divisor.

This turns intersections with $\widetilde{\mathcal R}$ into tangencies with
$\mathcal D_{\mathrm{red}}$.
Together with the positivity of the ramification divisor on the good model,
this is the bridge from algebraic geometry to counting: a section of large
height is forced to meet ramification, and such intersections become
tangencies to the branch divisor after projection to the abelian side.

\subsection{Toward the counting step}

The discussion above is still qualitative.  To finish the proof, one must
count the tangencies produced by Proposition~\ref{prop:tangency-branch}.  This
requires two further inputs.  The algebraic side, using the positivity of
$R_{\mathrm{can}}$ on the canonical model, and its pullback
$\widetilde R$ on the resolution, gives a lower bound for the tangency contribution.  The
analytic side, using Nevanlinna theory for the limiting linear entire curves,
gives an upper bound for the same contribution.  The contradiction is obtained
by comparing these two bounds.

There is one more technical point in passing from the qualitative
statement to the counting statement.  The fiber of $\mathcal D_{\mathrm{red}}$
over the limiting point $t_\infty\in T(\CC)$ may be non-reduced, and tangent
spaces need not specialize naively.  The following form of Gao's limiting
tangent-space lemma explains the proper closed subset which will appear in the
Nevanlinna estimate.

\begin{lemma}[Limiting tangent spaces]\label{lem:limiting-tangent-spaces}
Let $\mathcal D\to T$ be the family of branch divisors which occurs above, and
let $t_n\to t_\infty$.  After removing a proper closed subset from the reduced
limiting fiber $(\mathcal D_{t_\infty})_{\mathrm{red}}$, tangent directions to
the nearby reduced fibers $\mathcal D_{t_n,\mathrm{red}}$ specialize to tangent
directions of $(\mathcal D_{t_\infty})_{\mathrm{red}}$.
\end{lemma}

Thus tangencies to nearby branch divisors can be passed to the
limiting reduced divisor away from this proper closed subset.  The contribution
of the excluded subset will later be controlled by the codimension-two part of
the Second Main Theorem input.  This is the role of
\cite[Lemma 3.7]{GaoCounting}.

\section{Nevanlinna theory on expanding discs}\label{sec:nevanlinna-expanding-discs}

Before applying the counting argument, we isolate the part of Nevanlinna
theory that is used in the proof.  This is worthwhile for two reasons.  First,
the notation is not standard for every reader.  Second, the application here
is not quite the classical situation of one fixed entire curve: the curves
come from reparametrized sections on discs whose radii tend to infinity, and
only in the limit do we obtain an entire curve.  We use the notation of
\cite[\S 3.1]{GaoCounting}, but keep the discussion at the level needed for
the survey.

\subsection{Counting functions}

Let $Y$ be a projective variety over $\CC$, let
$\Delta_R=\{z\in\CC:|z|<R\}$, and let
\[
  \phi:\Delta_R\longrightarrow Y
\]
be a holomorphic map.  We first define the counting function for a closed
subscheme $Z\subset Y$, assuming that $\phi(\Delta_R)$ is not contained in
$Z$.

Choose effective Cartier divisors $E_1,\ldots,E_m$ on $Y$ such that
\[
  I_Z=I_{E_1}+\cdots+I_{E_m}.
\]
Equivalently, $Z=E_1\cap\cdots\cap E_m$ as a closed subscheme.  For
$z\in\Delta_R$, define
\[
  \operatorname{ord}_z\phi^*Z
  =
  \min_i \operatorname{ord}_z\phi^*E_i.
\]
This number is independent of the chosen Cartier divisors.  For $1<r<R$, the
counting function is
\[
  N(r,\phi,Z)
  =
  \int_1^r
  \left(
    \sum_{z\in\Delta_t}\operatorname{ord}_z\phi^*Z
  \right)\frac{dt}{t}.
\]
More generally, for an integer $k\ge 1$, the truncated counting function is
\[
  N^{(k)}(r,\phi,Z)
  =
  \int_1^r
  \left(
    \sum_{z\in\Delta_t}
      \min\{k,\operatorname{ord}_z\phi^*Z\}
  \right)\frac{dt}{t}.
\]
The case $k=1$ is the most important one below.  If $Z=E$ is a reduced
divisor, then
\[
  N(r,\phi,E)-N^{(1)}(r,\phi,E)
\]
measures the excess multiplicity of the intersections with $E$.  In geometric
terms, this is the Nevanlinna-theoretic quantity which records tangencies.

\subsection{Characteristic and proximity functions}

Let $L$ be a line bundle on $Y$ equipped with a smooth Hermitian metric
$\|\cdot\|$.  The characteristic function of $\phi$ with respect to $L$ is
\[
  T(r,\phi,L)
  =
  \int_1^r
  \left(
    \int_{\Delta_t}\phi^*c_1(L,\|\cdot\|)
  \right)\frac{dt}{t}.
\]
It measures the area growth of the holomorphic map with respect to the
curvature form of $L$.  If $E$ is an effective Cartier divisor, we also write
$T(r,\phi,E)$ for the characteristic function attached to $\mathcal O(E)$.

Now suppose that $E$ is an effective Cartier divisor defined by a section
$s_E\in H^0(Y,\mathcal O(E))$, and that $\phi(\Delta_R)$ is not contained in
$E$.  The proximity function is
\[
  m(r,\phi,E)
  =
  -
  \int_0^{2\pi}
  \log\|s_E(\phi(re^{i\theta}))\|
  \frac{d\theta}{2\pi}.
\]
Thus $m(r,\phi,E)$ measures how close the boundary circle
$\phi(\partial\Delta_r)$ is to $E$.  The counting function counts actual
intersections inside the disc; the proximity function measures near
intersections on the boundary.

\subsection{The First Main Theorem}

With the metric fixed, the First Main Theorem takes the following precise
form.  If $E$ is an effective Cartier divisor and
$\phi(\Delta_R)\not\subset E$, then for $1<r<R$,
\[
  T(r,\phi,E)
  =
  N(r,\phi,E)+m(r,\phi,E)-m(1,\phi,E).
\]
This identity is the analytic analogue of the elementary principle that the
total intersection growth with $E$ is split into actual intersections inside
the disc and boundary proximity to $E$.  In many presentations, different
choices of metrics change $T$ and $m$ only by bounded functions.  In the
argument below, the relevant metrics are fixed so that the error terms are
uniform for a varying sequence of maps.

\subsection{Limits of Nevanlinna functions}

We next record the point where the proof differs from the most standard
entire-curve setting.

\begin{lemma}[Nevanlinna functions under compact convergence]\label{lem:nevanlinna-limit}
Let $Y$ be a projective complex variety, and let
\[
  R_n\longrightarrow R_\infty\in(0,\infty]
\]
be a sequence of radii.  Let
\[
  \phi_n:\Delta_{R_n}\longrightarrow Y,\qquad 1\le n\le\infty,
\]
be holomorphic maps, where $\Delta_\infty=\CC$.  Suppose that
$\phi_n\to\phi_\infty$ uniformly on compact subsets of
$\Delta_{R_\infty}$.  Then for every fixed $r<R_\infty$,
\[
  T(r,\phi_n,L)\longrightarrow T(r,\phi_\infty,L).
\]
If $E$ is an effective divisor and
$\phi_\infty(\partial\Delta_r)\cap E=\emptyset$, then
\[
  N(r,\phi_n,E)\longrightarrow N(r,\phi_\infty,E),
  \qquad
  m(r,\phi_n,E)\longrightarrow m(r,\phi_\infty,E).
\]
For a closed subscheme $Z$ of arbitrary codimension, the useful replacement is
the one-sided estimate
\[
  N(r,\phi_n,Z)\le 2N(r,\phi_\infty,Z)
\]
for all sufficiently large $n$, depending on $r$.
\end{lemma}

This is the form of \cite[Lemma 3.5]{GaoCounting} that will be used below.  It
allows one to prove estimates on the limiting entire curve and then transfer
them back to the approximating maps on large finite discs.

\subsection{Yamanoi's Second Main Theorem}

The second input is the Second Main Theorem for abelian varieties, together
with its codimension estimate.  We use the following form.

\begin{theorem}[Second Main Theorem input]\label{thm:smt-input}
Let $A_0$ be a complex abelian variety, let $L$ be an ample line bundle on
$A_0$, and let $\psi:\CC\to A_0$ be a holomorphic curve with Zariski dense
image.
\begin{enumerate}
\item If $E$ is a reduced effective divisor on $A_0$, then for every
$\epsilon>0$,
\[
  T(r,\psi,E)
  \le
  N^{(1)}(r,\psi,E)+\epsilon T(r,\psi,L)
\]
outside a subset of $\mathbb R_{>0}$ of finite Lebesgue measure.
\item If $Z\subset A_0$ is a closed subscheme whose support has codimension at
least two, then for every $\epsilon>0$,
\[
  N(r,\psi,Z)\le \epsilon T(r,\psi,L)
\]
outside a subset of $\mathbb R_{>0}$ of finite Lebesgue measure.
\end{enumerate}
\end{theorem}

This is the form recalled in \cite[Theorem 3.6]{GaoCounting}, based on
Yamanoi's work.  Together with the First Main Theorem, the first part says
that multiple intersections with a reduced divisor are negligible compared
with the main characteristic growth.  We will use the theorem only through the
consequence below.

\subsection{A usable upper estimate}

We now package the Nevanlinna input in the form needed by the counting
argument.  Let $\mathcal D_{\mathrm{red}}\subset\mathcal A\times_\CC T$ be the
reduced branch divisor from the good model package.  Suppose that
\[
  \varphi_n:\Delta_{r_n}\longrightarrow \mathcal A\times_\CC T,
  \qquad r_n\to\infty,
\]
converges uniformly on compact subsets to a linear entire curve
\[
  \varphi_\infty:\CC\longrightarrow \mathcal A_b\times\{t_\infty\}.
\]
Write
\[
  D_\infty=(\mathcal D_{b,t_\infty})_{\mathrm{red}}
  \subset \mathcal A_b\times\{t_\infty\}.
\]
The Second Main Theorem will be applied on the Zariski closure
$V_\infty$ of $\varphi_\infty(\CC)$, rather than on the whole abelian fiber.
After translating $V_\infty$ to the abelian subvariety through the origin, the
limiting curve becomes Zariski dense in that abelian variety, so the hypotheses
of Theorem~\ref{thm:smt-input} apply to the restrictions of the divisors and
closed subschemes to $V_\infty$.

\begin{proposition}[Nevanlinna upper estimate for expanding discs]\label{prop:nevanlinna-upper}
In the set-up above, assume the following additional conditions.
\begin{enumerate}
\item Let $V_\infty$ be the Zariski closure of
$\varphi_\infty(\CC)$ in $\mathcal A_b\times\{t_\infty\}$.  Then
$V_\infty$ is a translate of an abelian subvariety, and we regard
$\varphi_\infty$ as a Zariski dense entire curve in $V_\infty$.
\item The divisor $D_\infty$ does not contain $V_\infty$.  In particular,
the counting functions
\[
  N(r,\varphi_\infty,D_\infty)
  \quad\text{and}\quad
  N^{(1)}(r,\varphi_\infty,D_\infty)
\]
are defined.  Moreover, after discarding finitely many terms, the maps
$\varphi_n$ are not contained in $\mathcal D_{\mathrm{red}}$, so the
counting functions involving $\varphi_n$ below are also defined.
\end{enumerate}
Then for every $\epsilon>0$, outside a subset of $\mathbb R_{>0}$ of finite
Lebesgue measure, one has
\[
  N(r,\varphi_n,\mathcal D_{\mathrm{red}})
  -N^{(1)}(r,\varphi_n,\mathcal D_{\mathrm{red}})
  \le
  \epsilon\,T(r,\varphi_\infty,D_\infty)
\]
for all sufficiently large $n$, depending on $r$.
\end{proposition}

The point of the proposition is that tangencies of the approximating discs to
the reduced branch divisor have asymptotic proportion zero.  The proof uses
Theorem~\ref{thm:smt-input} on the limiting linear entire curve, transfers the
estimate to the maps $\varphi_n$ by Lemma~\ref{lem:nevanlinna-limit}, and
removes the proper closed subset described in
Lemma~\ref{lem:limiting-tangent-spaces}, where tangent spaces of
nearby reduced fibers may fail to specialize well.  The contribution of this closed subset is
controlled by Theorem~\ref{thm:smt-input}(2), applied to a closed subset whose
support has codimension at least two.  In the
application below, the two assumptions are ensured by the very general
choice of the fiber, the mobility of limiting linear entire curves, and the
genericity of the original sequence.

\section{The counting argument and the proof of the main theorem}\label{sec:counting-proof}

We now assemble the preceding ingredients and finish the proof in the
finite-cover case.  Let us first recall what has been reduced.  The main
result was stated in several related forms.  The geometric form is
Theorem~\ref{thm:finite-cover-geometric}; for the proof it is enough to
establish the generic-sequence form, Theorem~\ref{thm:finite-cover-generic}.
In the curve-function-field setting $K=\CC(B)$, where $B$ is a smooth
projective complex curve, the assertion is that if $X/K$ is a projective
variety of general type admitting a finite morphism to an abelian variety,
then every generic sequence in $X(K)$ has bounded height.

Section~\ref{sec:gao-counting-setup} reduced this statement to
Proposition~\ref{prop:gao-generic-reduction}.  Thus, after the trace and
quotient reductions, and after replacing the original data by the good model
package of Proposition~\ref{prop:good-model-package}, we are in the following
situation.  There is a finite surjective morphism
\[
  f:X\longrightarrow A\times_KT_K,
\]
where $A/K$ is an abelian variety and $T$ is a projective variety over $\CC$.
We are given a generic sequence $\{x_n\}\subset X(K)$ such that
\[
  f(x_n)=(s_n,t_n),
  \qquad s_n\in A(K),\quad t_n\in T(\CC),
\]
and we must prove that the heights of $x_n$ are bounded.  This is the content
of Proposition~\ref{prop:gao-generic-reduction}.

We argue by contradiction.  Suppose that the heights of $x_n$ tend to
infinity.  The construction of limiting linear entire curves, together with
their mobility, then produces, after passing to a subsequence, choosing a very
general point $b\in B(\CC)$, and reparametrizing near $b$, a sequence of
holomorphic maps
\[
  \varphi_n:\Delta_{r_n}\longrightarrow \mathcal A\times_\CC T,
  \qquad r_n\to\infty,
\]
which converges uniformly on compact subsets to a limiting linear entire curve
\[
  \varphi_\infty:\CC\longrightarrow \mathcal A_b\times\{t_\infty\}.
\]
The contradiction will come from counting how often the maps $\varphi_n$ are
tangent to the reduced branch divisor $\mathcal D_{\mathrm{red}}$.

\subsection{The expanding-disc set-up}

Since $T$ is projective, the constant points $t_n$ have a convergent
subsequence; we have denoted its limit by $t_\infty$.  The maps
$\varphi_n$ above are the reparametrized images of the sections
$(s_n,t_n)$, and the limiting map lies over the limiting constant point
$t_\infty$.

The divisor on the limiting fiber is
\[
  D_\infty=(\mathcal D_{b,t_\infty})_{\mathrm{red}}
  \subset \mathcal A_b\times\{t_\infty\}.
\]
The growth term against which all estimates are compared is
$T(r,\varphi_\infty,D_\infty)$.  The quantity to be estimated is
\[
  N(r,\varphi_n,\mathcal D_{\mathrm{red}})
  -N^{(1)}(r,\varphi_n,\mathcal D_{\mathrm{red}}),
\]
which measures the tangency contribution of $\varphi_n$ along the reduced
branch divisor.

\subsection{The Nevanlinna upper bound}

The first estimate has already been packaged in
Proposition~\ref{prop:nevanlinna-upper}.  It says that the tangency
contribution of the approximating discs along the reduced branch divisor is
negligible:
\[
  N(r,\varphi_n,\mathcal D_{\mathrm{red}})
  -N^{(1)}(r,\varphi_n,\mathcal D_{\mathrm{red}})
  \le
  \epsilon\,T(r,\varphi_\infty,D_\infty)
\]
outside an exceptional set of radii of finite measure, for all large $n$.
Thus the analytic side gives an upper bound of proportion zero.

\subsection{The ramification lower bound}

The opposite estimate comes from the algebraic geometry of
Section~\ref{sec:gao-counting-setup}.  The positivity of
$R_{\mathrm{can}}$, pulled back to $\widetilde R$ on the resolution,
forces the lifted sections to meet ramification often.
Proposition~\ref{prop:tangency-branch} then converts these ramification
intersections into tangencies with $\mathcal D_{\mathrm{red}}$.

\begin{proposition}[Ramification lower bound]\label{prop:ramification-lower}
There is a constant $c>0$ and a subset of $\mathbb R_{>0}$ of finite Lebesgue
measure such that, outside this subset,
\[
  N(r,\varphi_n,\mathcal D_{\mathrm{red}})
  -N^{(1)}(r,\varphi_n,\mathcal D_{\mathrm{red}})
  \ge
  c\,T(r,\varphi_\infty,D_\infty)
\]
for all sufficiently large $n$, depending on $r$.
\end{proposition}

This is the content of \cite[Lemma 3.18]{GaoCounting}.  Its proof has three
ingredients.  First, the relative positivity of $R_{\mathrm{can}}$ on the
canonical model gives, after pullback to $\widetilde X$, a positive
lower bound for intersections with the lifted
maps.  Second, the tangency statement converts the relevant part of these
intersections into the excess term
$N-N^{(1)}$ for the reduced branch divisor.  Third, the good-model reductions
control the remaining ramification components: the relevant branch divisor is
horizontal over the constant factor, while the exceptional higher-codimension
terms are negligible for the limiting entire curve by
Theorem~\ref{thm:smt-input}(2).  Thus the algebraic side says that tangencies
occur with positive proportion.

\subsection{The contradiction}

We now finish the proof of the key boundedness statement, and hence of
Theorem~\ref{thm:finite-cover-generic}.  Choose
\[
  0<\epsilon<c,
\]
where $c$ is the constant in Proposition~\ref{prop:ramification-lower}.  The
exceptional sets of radii in Propositions~\ref{prop:nevanlinna-upper} and
\ref{prop:ramification-lower} have finite Lebesgue measure, so we may choose
$r$ outside their union and large enough that
\[
  T(r,\varphi_\infty,D_\infty)>0.
\]
For this fixed $r$ and all sufficiently large $n$, the two propositions give
\[
  c\,T(r,\varphi_\infty,D_\infty)
  \le
  N(r,\varphi_n,\mathcal D_{\mathrm{red}})
  -N^{(1)}(r,\varphi_n,\mathcal D_{\mathrm{red}})
  \le
  \epsilon\,T(r,\varphi_\infty,D_\infty),
\]
which is impossible because $\epsilon<c$ and
$T(r,\varphi_\infty,D_\infty)>0$.  This contradiction shows that the assumed
generic sequence of unbounded height cannot exist, completing the proof of
Proposition~\ref{prop:gao-generic-reduction}.

\section{Further directions}\label{sec:further-directions}

The proof explained in this survey suggests a broader question: how far can one
make Vojta's dictionary concrete in Diophantine problems?  In the
function-field setting considered here, the guiding philosophy is that a
height-unbounded sequence of rational points should, after a suitable limiting
process, produce analytic objects whose existence is constrained by
hyperbolicity or Nevanlinna theory.  In the
abelian-variety setting considered above, this philosophy is implemented using
three rather specific ingredients: the group structure of rational points, the
partial-height construction which avoids bad fibers, and a ramification
counting argument.  These ingredients should not be viewed as the only possible
implementation.  In a more directly hyperbolic setting, for instance, one may
not need such refined counting estimates; in a more general setting, one may
instead need substitutes for one or more of the abelian-variety inputs, or even
a different way to realize the same dictionary.

A first direction is to look beyond abelian varieties, especially toward
classes of varieties for which the analytic Green--Griffiths--Lang picture is
already known.  General hypersurfaces of sufficiently large degree in
projective space provide a basic testing ground.  Their hyperbolicity was
established in important work of Siu \cite{SiuHyperbolicity}; later
developments, including work of Brotbek and of Cadorel, give stronger and
more flexible forms of the hyperbolicity and Green--Griffiths--Lang
statements \cite{BrotbekHyperbolicity,CadorelHypersurfaces}.  In such
examples the analytic obstruction to entire curves is well understood.  The
missing ingredient, from the Diophantine side, is a replacement for the
abelian-group structure used above: it produces limiting linear entire curves
with enough mobility, and it is also what makes the partial-height localization
non-degenerate.

A second direction stays closer to abelian varieties but changes the
Diophantine problem.  Since rational points on abelian varieties are highly
structured, one may ask whether partial heights and limiting analytic objects
can also be useful for other function-field problems involving abelian
varieties, such as Campana points, gcd-type problems, or abc-type inequalities;
see, for example, the perspective on Campana points in
\cite{AbramovichVarillyCampana}.  These questions are not the same as the
projective Bombieri--Lang problem, but they share a similar feature: one wants
to control intersections with boundary divisors or exceptional loci.  The
method surveyed here suggests that, when the points involved still have
enough group-theoretic organization, height growth might again be converted
into analytic limits and then into counting estimates.

A third possible arena is arithmetic dynamics over function fields.  In this
setting the role of a finitely generated subgroup may be played by the orbit
of an endomorphism.  Such orbits are not groups, but they are still highly
organized sequences of rational points.  This makes it natural to ask whether
height-growing orbits can produce limiting analytic objects analogous to the
linear entire curves used above.  If so, one might hope to apply related ideas
to classical problems such as the function-field dynamical Mordell--Lang
problem, where one studies the intersection of an orbit with a subvariety; see
\cite{BellGhiocaTuckerDML,XieAroundDML} for the general dynamical
Mordell--Lang problem and related developments.

There is also a question about what happens near bad fibers.  In the proof
discussed above, partial heights are used to localize the construction on good
analytic discs of the base curve.  This avoids degeneration and leads to
complex entire curves on smooth fibers.  If one allows the area to concentrate
near bad fibers instead, the limiting object may no longer be an ordinary
complex entire curve.  At present, however, non-archimedean hyperbolicity
theory is still much less developed than its complex counterpart, and even the
correct analogue of a complex entire curve is not completely clear.  It is
nevertheless tempting to ask whether some non-archimedean or Berkovich
limiting object should appear in such a situation, and whether it could carry
useful Diophantine information.

Finally, one may ask whether any part of this picture has a number-field
analogue.  Over a function field, rational points are sections over a curve,
and height growth can be converted into area growth.  This geometric feature
is absent over number fields.  Even for finite covers of abelian varieties,
the corresponding number-field case of Bombieri--Lang remains wide open.
Thus any number-field analogue of the method would require a substitute for
the passage from sections to limiting analytic curves.  Understanding what
this substitute should be is itself part of the broader problem.

\section*{Acknowledgements}

This article is based on a talk given at the International Congress of Chinese
Mathematicians.  I thank the organizers of the ICCM for the invitation.  I am
supported by the National Natural Science Foundation of China (Grant
No.~12271007).  I thank Guoquan Gao and Junjiro Noguchi for useful comments on the first
version of this article.  I thank the AI assistant Xiaozhua for assistance during the
writing process, and Liang Xiao and Shuai Chen for providing the AI-agent
environment and for resolving related technical issues.

\bibliographystyle{alpha}
\bibliography{gbl_refs}

\end{document}